\documentclass{amsart}
\usepackage{amsmath}
\usepackage{amssymb,amsthm}
\usepackage{graphicx}
\usepackage{tikz}
\usepackage{mathrsfs}
\numberwithin{equation}{section}
\newtheorem{theorem}{Theorem}[section]
\newtheorem{remark}[theorem]{Remark}

\newtheorem{lemma}[theorem]{Lemma}
\newtheorem{corollary}[theorem]{Corollary}
\newtheorem{definition}[theorem]{Definition}
\newtheorem{proposition}[theorem]{Proposition}

\title[solutions to the Hartree-Fock equation]{Structures of sets of solutions to the Hartree-Fock equation}
\author{Sohei Ashida}
\begin{document}
\maketitle

\begin{abstract}
The Hartree-Fock equation which is the Euler-Lagrange equation corresponding to the Hartree-Fock energy functional is used in many-electron problems. Since the Hartree-Fock equation is a system of nonlinear eigenvalue problems, the study of structures of sets of all solutions needs new methods different from that for the set of eigenfunctions of linear operators. In this paper we prove that the sets of all solutions to the Hartree-Fock equation associated with critical values of the Hartree-Fock energy functional less than the first energy threshold are unions of a finite number of compact connected real-analytic spaces. The result would also be a basis for the study of approximation methods to solve the equation.
\end{abstract}

\section{Introduction}
In this paper we study sets of solutions to the Hartree-Fock equation which is used to seek approximated eigenfunctions of the electronic Hamiltonian by the variational method and afford functions used in various other approximation methods. Fix the number of electrons $N\in \mathbb N$ satisfying $N\geq 2$, number of nuclei $n\in \mathbb N$, nuclear charges $Z_j>0,\ j=1,\dots,n$, and the positions of the nuclei $\bar x_j\in \mathbb R^3,\ j=1,\dots,n$. Let $\varphi_i\in H^2(\mathbb R^3),\ i=1,\dots, N$ and set $\Phi:={}^t(\varphi_1,\dots,\varphi_N)$. The function $\varphi_i$ is complex-valued, but everything in this paper is trivially adapted to spin-dependent functions with only notational changes. We define the Fock operator by $\mathcal F(\Phi):=h + R^{\Phi} - S^{\Phi},$
where $h :=-\Delta + V$ with $V(x):= -\sum_{j=1}^n\frac{Z_j}{\lvert x-\bar x_j\rvert}$, and operators depending on $\Phi$ are defined by $R^{\Phi}(x) :=\sum_{i=1}^N\int\lvert x-y\rvert^{-1}\lvert \varphi_i(y)\rvert^2 dy = \sum_{i=1}^NQ_{ii}^{\Phi}(x)$ and $S^{\Phi} := \sum_{i=1}^NS^{\Phi}_{ii}$
with
\begin{align*}
Q^{\Phi}_{ij}(x) &:=\int\lvert x-y\rvert^{-1}\varphi_j^*(y)\varphi_i(y)dy,\\
(S^{\Phi}_{ij}w)(x) &:= \left(\int \lvert x-y\rvert^{-1} \varphi_j^*(y) w(y) dy\right)\varphi_i(x).
\end{align*}
The Hartree-Fock equation is written as
\begin{equation}\label{myeq1.3}
\mathcal F(\Phi)\varphi_i = \epsilon_i \varphi_i,\ 1 \leq i \leq N,
\end{equation}
with constraints $\langle \varphi_i,\varphi_j\rangle=\delta_{ij},\ 1\leq i,j\leq N$, where $\epsilon_i\in \mathbb R,\ 1 \leq i \leq N$ are unknown constants.
We seek tuples of functions $(\varphi_1,\dots,\varphi_N)$ and constants $(\epsilon_1,\dots,\epsilon_N)$ satisfying \eqref{myeq1.3}. We call the tuple of the constants $(\epsilon_1,\dots,\epsilon_N)$ an orbital energy, if \eqref{myeq1.3} has a solution. 

The Hartree-Fock equation is the Euler-Lagrange equation corresponding to the Hartree-Fock energy functional $\mathcal E(\Phi) = \mathcal E_N(\Phi) := \langle \Psi,H\Psi\rangle$, where
\begin{equation*}
H:=-\sum_{i=1}^N\Delta_{x_i}+\sum_{i=1}^NV(x_i)+\sum_{1\leq i <j\leq N}\frac{1}{\lvert x_i-x_j\rvert},
\end{equation*}
is an electronic Hamiltonian acting on $L^2(\mathbb R^{3N})$, $\Phi={}^t(\varphi_1,\dots,\varphi_N) \in \bigoplus_{i=1}^N H^1(\mathbb R^3)$ with constraints $\langle \varphi_i,\varphi_j\rangle=\delta_{ij},\ 1\leq i,j\leq N$, and $\Psi$ is the Slater determinant
\begin{equation*}
\Psi(x_1,\dots,x_N):=(N!)^{-1/2}\sum_{\tau\in\mathbf S_N}(\mathrm{sgn}\, \tau)\varphi_1(x_{\tau(1)})\cdots\varphi_N(x_{\tau(N)}).
\end{equation*}
Here $\mathbf S_N$ is the symmetric group and $\mathrm{sgn}\, \tau$ is the signature of $\tau$. The functional $\mathcal E(\Phi)$ can be written explicitly as
\begin{equation*}\label{myeq1.2}
\mathcal E(\Phi)=\sum_{i=1}^N\langle \varphi_i ,h\varphi_i\rangle + \frac{1}{2}\int\int\rho(x)\frac{1}{\lvert x-y\rvert}\rho(y)dxdy - \frac{1}{2}\int\int\frac{1}{\lvert x-y\rvert}\lvert\rho(x,y)\rvert^2 dx dy,
\end{equation*}
where $\rho(x) := \sum_{i=1}^N\lvert \varphi_i(x)\rvert^2$ is the density, and $\rho(x,y) := \sum_{i=1}^N \varphi_i(x) \varphi_i^*(y)$ is the density matrix. The relation between the Hartree-Fock equation and $\mathcal E(\Phi)$ is as follows. Solutions to the Hartree-Fock equation are critical points of $\mathcal E(\Phi)$, and $E\in \mathbb R$ is a critical value of $\mathcal E(\Phi)$  if and only if there exists a solution $\Phi$ to the Hartree-Fock equation such that $\mathcal E(\Phi)=E$. Although $\mathcal E(\Phi)$ is defined on $\bigoplus_{i=1}^NH^1(\mathbb R^3)$, by the standard regularity results the critical points satisfying the Hartree-Fock equation belong to $\bigoplus_{i=1}^NH^2(\mathbb R^3)$ (see e.g. \cite{Li}). The Hartree-Fock equation was introduced independently by Fock \cite{Fo} and Slater \cite{Sl}, after Hartree \cite{Ha} introduced the Hartree equation ignoring the antisymmetry with respect to exchanges of variables. 

Lieb-Simon \cite{LS} proved that if $N<\sum_{j=1}^nZ_j+1$, there exists a solution to the Hartree-Fock equation which minimizes the functional. Lions \cite{Li} proved that, if $N\leq \sum_{j=1}^nZ_j$, there exists a sequence of solutions to the Hartree-Fock equation with nonpositive orbital energies such that the corresponding critical values are converging to $0$. Lewin \cite{Le} showed that there exists a sequence of critical values of the Hartree-Fock functional less than the first energy threshold $J(N-1)$ and converging to $J(N-1)$ under the same assumption as in \cite{LS}, where
\begin{equation}\label{myeq1.4}
J(N-1):=\inf \{\mathcal E_{N-1}(\Phi) :\Phi \in \bigoplus_{i=1}^{N-1}H^1(\mathbb R^3),\ \langle \varphi_i,\varphi_j\rangle =\delta_{ij},\ 1\leq i,j\leq N-1\}.
\end{equation}
Ashida \cite{As} proved that for any $\epsilon>0$ the set of all critical values less than $J(N-1)-\epsilon$ is finite for arbitrary $N$. If we ignore the antisymmetry, the Hartree equation instead of the Hartree-Fock equation is obtained from the variational problem for electronic Hamiltonians. For the existence of the solutions to the Hartree equation see e.g. \cite{Re,Wo,St,LS,Li}.

Even for small $N$ and $n$ the Hartree-Fock equation can not be solved exactly. In approximation methods to obtain solutions to the Hartree-Fock equation such as self-consistent-field (SCF) method, we consider convergence of a sequence of functions obtained successively in the methods. However, without knowing the structure of the set of all solutions to the Hartree-Fock equation, arguments about the convergence would necessarily be unsatisfactory from both theoretical and practical point of views and the results would be rather weak statements (e.g. see the statement in \cite[Theorem 7]{CB} and the comments below the theorem), because the limit points of the converging sequences are unknown objects. The structure of the set of critical points is more complicated than that of critical values. We can illustrate this fact from the simple example $f:\mathbb R^2\to\mathbb R$ defined by $f(x,y)=x^2$ whose critical value is $0$ and critical points are $\{(0,y): y\in \mathbb R\}$. We prove in this paper that the set $A(E)$ of all solutions to the Hartree-Fock equation associated with a critical value $E<J(N-1)$ has a structure as a union of a finite number of compact connected real-analytic spaces without isolated points. Let us compare this result with eigenvalue problems of linear operators. In the case of a linear operator $H$ on a Hilbert space $X$, the critical values of the functional $\langle u,Hu\rangle$ for $u\in X$ with the constraint $\lVert u\rVert =1$ are eigenvalues and the corresponding critical points are eigenfunctions. If the eigenvalue is simple the corresponding set of normalized eigenfunctions is $\{cu: c\in \mathbb C,\ \lvert c\rvert =1\}$, where $u$ is a normalized eigenfunction. This set is a circle in the complex plane. When the multiplicity of the eigenvalue is finite, the set of all normalized eigenfunctions is a finite dimensional sphere. Thus the present result is analogous to the finite multiplicity case of linear operators. For any $\epsilon>0$ we also obtain a similar result for the set $B(\epsilon)$ of all solutions associated with orbital energies $(\epsilon_1,\dots,\epsilon_n)$ satisfying $\epsilon_i<-\epsilon,\ 1\leq i\leq N$.

For the proof of the main result, we consider an auxiliary functional $f$ closely related to $\mathcal E$. The set of all solutions to the Hartree-Fock equation associated with a critical value $E<J(N-1)$ is regarded as a subset of the set of all critical points of $f$ associated with the same critical value $E$. At each critical point of $f$ we show that there exists a neighborhood in the whole space and a subset of the neighborhood homeomorphic to an open set in a finite-dimensional Euclidean space including the set of all critical points in the neighborhood. The subset can be regarded as a manifold. This is achieved by showing that the Fr\'echet second derivative of $f$ is a Fredholm operator and using the implicit function theorem in Banach spaces. The first derivative is defined using a bilinear form as in Fu\v cik-Ne\v cas-Sou\v cek-Sou\v cek \cite{FNSS} rather than by the usual definition of the Fr\'echet derivative. The merit of this method in our proof is to reduce the complexity due to inner products which takes complex conjugate of the functions. The Fredholm property of the second derivative is proved showing that the derivative is a sum of an isomorphism and a compact operator. The importance of this decomposition of the derivative was discovered by \cite{FNSS} in order to prove that a critical value of the functional is isolated in the set of all critical values near a critical point (see in particular \cite[Theorem 4.1]{FNSS}) and apply the result to an eigenvalue problem for nonlinear elliptic equations in a bounded domain with a Dirichlet boundary condition.  Next, we show that the set of all solutions to the Hartree-Fock equation associated with $E$ in the manifold is the set of all zeros of some real-analytic functions. Finally, we show that in the intersection of the manifolds the mapping between two real-analytic subsets is a mapping of ringed spaces, which can be achieved by a natural extension of the mapping to some real-analytic mapping defined on the manifold. The real-analyticity of the mapping follows from the real-analyticity of the mapping in the implicit function theorem which is a consequence of the real-analyticity of $f$.

As mentioned above, a standard approximation method to obtain a solution to the Hartree-Fock equation is SCF method. Convergence of SCF method is important for practical calculations to solve the Hartree-Fock equation. SCF method is rarely studied from a rigorous standpoint. One of a few results is the work by \cite{CB} proving that under a certain assumption there exists a subsequence of the functions either converging to a critical point of the Hartree-Fock functional or oscillating between two states. We would need to study the relation between these limit points and the set of solutions and how the structure of the set of solutions affects the convergence. The present result about the structure of sets of solutions would be a basis for such a study of the convergence problems.

The content of this paper is as follows. In Sect. \ref{secondsec} we state our main results. In Sect. \ref{thirdsec} we introduce the notion of real-analytic operators in Banach spaces and give lemmas needed in the proof of the main results as preliminaries. In Sect. \ref{fourthsec} we prove the main theorem. In Appendix we prove a lemma in the preliminary.

\section{Main results}\label{secondsec}
Before we state the results precisely, we shall make clear what is meant by a critical value of the Hartree-Fock functional $\mathcal E$ and the relation between $\mathcal E$ and the Hartree-Fock equation. If there is no additional constraints, a critical point of a functional is defined as a point at which the Fr\'echet derivative of the functional vanishes. Considering the functional $\mathcal E(\Phi)$ with constraints $\langle \varphi_i,\varphi_j\rangle=\delta_{ij}$ is equivalent to restricting the functional to the subset $\mathcal N :=\{\Phi\in\bigoplus_{i=1}^NH^1(\mathbb R^3): \langle \varphi_i,\varphi_j\rangle=\delta_{ij}\}$. In the following definition of the critical point we may suppose more generally $\mathcal N$ is a subset of a locally convex space $X$ and $\mathcal E(\Phi):X\to \mathbb R$. Let $\Phi\in\mathcal N$, and $\mathcal C_{\Phi}$ be the set of all mappings $c:(-1,1)\to \bigoplus_{i=1}^NH^1(\mathbb R^3)$ such that $c(t)\in \mathcal N$ for any $t\in (-1,1)$, $c(0)=\Phi$ and $c'(0)$ exists. Then $\Phi$ is a critical point of $\mathcal E(\Phi)$ if and only if $\frac{d\mathcal E(c(t))}{dt}\big\vert_{t=0}=0$ for any $c\in \mathcal C_{\Phi}$ (cf \cite[Definition 43.20]{Ze2}).  A real number $E$ is called a critical value if there exists a critical point $\Phi\in \mathcal N$ such that $\mathcal E(\Phi)=E$.

In order to seek critical points, we have the method of Lagrange multiplier as in finite-dimensional cases (cf \cite[Proposition 43.21]{Ze2}). Let us define a functional $g_{ij}(\Phi):=\langle\varphi_i,\varphi_j\rangle-\delta_{ij}$ (note that for $i\neq j$, $g_{ij}$ is complex-valued and $g_{ij}=g_{ji}^*$ ). Then by the method of Lagrange multiplier we can see that $\Phi\in\mathcal N$ is a critical point of $\mathcal E(\Phi)$, if and only if there exists an Hermitian matrix $(\epsilon_{ij})$ such that $\Phi$ is a critical point of the functional
\begin{align*}
&\mathcal E(\Phi)-\sum_{1\leq i< j\leq N}2(\mathrm{Re}\, \epsilon_{ij}\mathrm{Re}\, g_{ij}(\Phi)-\mathrm{Im}\, \epsilon_{ij}\mathrm{Im}\, g_{ij}(\Phi))-\sum_{i=1}^N\epsilon_{ii}g_{ii}(\Phi)\\
&\quad=\mathcal E(\Phi)-\sum_{1\leq i,j\leq N}\epsilon_{ij}g_{ij}(\Phi),
\end{align*}
without any restriction. Therefore, by a direct calculation of the Fr\'echet derivative we can see that a sufficient and necessary condition for $\Phi$ to be a critical point of $\mathcal E(\Phi)$ is that $\Phi$ satisfies the following equation
$$\mathcal F(\Phi)\varphi_i = \sum_{j=1}^N\epsilon_{ij} \varphi_j,\ 1 \leq i \leq N,$$
for some Hermitian matrix $(\epsilon_{ij})$. Since $(\epsilon_{ij})$ is Hermitian, it is diagonalized by a unitary $N\times N$ matrix $(a_{ij})$. Thus after a unitary change $\varphi_i^{\mathrm{New}}=\sum a_{ij}\varphi_j$,  the Hartree-Fock equation \eqref{myeq1.3} is satisfied by $(\varphi_1^{\mathrm{New}},\dots,\varphi_N^{\mathrm{New}}),$
and some real numbers $(\epsilon_1,\dots,\epsilon_N)$. Therefore, a real number $E$ is a critical value of $\mathcal E(\Phi)$ if and only if there exists a solution $\Phi$ to the Hartree-Fock equation such that $\mathcal E(\Phi)=E$.

To express the structure of the set of solutions to the Hartree-Fock equation, we introduce a few standard notions. We call a Hausdorff space $M$ a real-analytic manifold, if for some $n\in \mathbb N$, $M$ is locally homeomorphic to $n$-dimensional  Euclidean space and its coordinate transformations are real-analytic. We call a subset $A$ of a real-analytic manifold $M$ a real-analytic subset, if there exists an open cover $\{U_i\}$ of $M$ such that for any $i$ there exist real-analytic functions $f_i^1,\dots, f_i^s$ on $U_i$ satisfying $A\cap U_i=\{z\in U_i : f_i^1(z)=\dotsm=f_i^s(z)=0\}$. Although the set of solutions to the Hartree-Fock equation has a structure of real-analytic subset in a neighborhood of each point, in order to express the whole structure we need the notion of real-analytic space (see e.g. \cite{GR, GMT}). Correspondence of the structure of real-analytic subsets in an intersection of neighborhoods is judged by the correspondence of the sheaf of germs of real-analytic functions on the real-analytic subsets. A ringed space is a pair $(A,{}_A\mathscr O)$ of a Hausdorff space $A$ and a sheaf ${}_A\mathscr O$ of subrings of the sheaf of germs of continuous real-valued functions on $A$. A ringed space $(A,{}_A\mathscr O)$ is called a real-analytic space if every $x\in A$ has a neighborhood $W\subset A$ such that $(W,{}_A\mathscr O|_W)$ is isomorphic to a ringed space $(A',{}_{A'}\mathscr O)$, where $A'$ is a real-analytic subset of a domain $V$ in $\mathbb R^{\mu},\ \mu\in\mathbb N$ and ${}_{A'}\mathscr O=(\mathscr O/\mathscr I)|_{A'}$ (This kind of space is called real analytic variety in \cite{GMT} and reduced real analytic space in \cite{ABF}). Here $\mathscr O$ is the sheaf of germs of real-analytic functions over $V$ and $\mathscr I$ is the sheaf of germs of real-analytic functions that vanish on $A'$. The reason to take the quotient $\mathscr O/\mathscr I$ is to avoid the redundancy due to functions which coincide on $A'$ but different on the complementary set $V\setminus A'$, when they are restricted to $A'$. Real-analytic space is a generalization of real-analytic manifold which is allowed to have singular points. An example of the singular point is $(x,y)=(0,0)$ in $\{(x,y)\in\mathbb R^2: x^2-y^2=0\}$. We often denote an real-analytic space $(A,{}_A\mathscr O)$ simply by $A$.

Let $A(E)$ be the set of all solutions $\Phi$ to the Hartree-Fock equation \eqref{myeq1.3} associated with a critical value $E$ of the Hartree-Fock functional i.e. solutions such that $\mathcal E(\Phi)=E$. We consider the relative topology induced by the topology of $\bigoplus_{i=1}^NH^2(\mathbb R^3)$ as the topology of $A(E)$. Our main results are the followings.
\begin{theorem}\label{asub}
For any critical value $E<J(N-1)$ of the Hartree-Fock functional, $A(E)$ is a union $\bigcup_{k=1}^mA_k(E)$ of a finite number of compact connected real-analytic spaces $A_k(E),\ k=1,\dots, m$, and $A(E)$ does not have isolated points.
\end{theorem}
For $\epsilon>0$, let $B(\epsilon)$ be the set of all solutions to the Hartree-Fock equation \eqref{myeq1.3} associated with orbital energies $(\epsilon_1,\dots,\epsilon_N)$ satisfying $\epsilon_i<-\epsilon,\ 1\leq i\leq N$.
\begin{theorem}\label{bsub}
For any $\epsilon>0$, $B(\epsilon)$ is a union $\bigcup_{k=1}^lB_k(\epsilon)$ of a finite number of connected real-analytic spaces $B_k(\epsilon),\ k=1,\dots l$, and $B(\epsilon)$ does not have isolated points.
\end{theorem}

We present an illustrative example of the result above. Since the Hartree-Fock equation itself can not be solved, we replace the potential $\sum_{1\leq i<j\leq N}\frac{1}{\lvert x_i-x_j\rvert}$ between electrons by $0$. Moreover, we assume $n=1$, $N=2$, $Z_1=1$ and $\bar x_1=0$. In this case the same results as above hold and we can solve the equation corresponding to the Hartree-Fock equation. The equation is written as
\begin{equation}\label{myeq2.1}
\begin{split}
&(-\Delta-\frac{1}{| x_i|})\varphi_i=\epsilon_i \varphi_i,\ i=1,2,\\
&\lVert \varphi_i\rVert=1,\ i=1,2,\ \langle \varphi_1,\varphi_2\rangle=0.
\end{split}
\end{equation}
If we ignore the constraint $\langle \varphi_1,\varphi_2\rangle=0$, the equations are just the eigenvalue problems for the Hydrogen atom and solved by separation of variables.
We have $J(N-1)=J(1)=-\frac{1}{4}$, and the critical values less than $J(N-1)$ are
$-\frac{1}{4}-\frac{1}{4j^2},\ j=2,3,\dots.$
The set of all solutions to the equation associated with the lowest critical value $-\frac{5}{16}$ is written as
\begin{equation}\label{myeq2.2}
A(-\frac{5}{16})=(\mathbb S^1\times\mathbb S^7)\sqcup(\mathbb S^7\times\mathbb S^1).
\end{equation}
Here $\mathbb S^{2d-1}$ is the unit sphere in the complex $d$-dimensional eigenspace of $-\Delta +\frac{1}{|x|}$. When the multiplicity of an eigenvalue $\epsilon$ is $d$, there exists an orthonormal basis $u_1,\dots, u_d$ of the eigenspace. Thus the set of all normalized solutions to $(-\Delta-\frac{1}{|x|})\varphi=\epsilon\varphi$ is the set $\{\sum_{i=1}^d c_iu_i: c_i\in \mathbb C,\ \sum_{i=1}^d|c_i|^2=1\}$. Setting $a_{2i-1}:=\mathrm{Re}\, c_i,\ a_{2i}:=\mathrm{Im}\, c_i$, this set is identified with the sphere $\mathbb S^{2d-1}=\{(a_1,a_2,\dots,a_{2d})\in \mathbb R^{2d}: \sum_{j=1}^{2d}a_j^2=1\}$. The spheres $\mathbb S^1$ and $\mathbb S^7$ above correspond to the multiplicities $1$ and $4$ of the first and the second eigenvalues of $-\Delta - \frac{1}{|x|}$ in ascending order. The two components in \eqref{myeq2.2} correspond to the two possibilities $\varphi_1\in\mathbb S^1,\ \varphi_2\in\mathbb S^7$ and $\varphi_1\in\mathbb S^7,\ \varphi_2\in\mathbb S^1$ (note that the constraint $\langle\varphi_1,\varphi_2\rangle=0$  prevents $\varphi_1,\varphi_2\in\mathbb S^1$ being a solution of \eqref{myeq2.1}).

If $\Phi$ is a solution to the Hartree-Fock equation, by the standard regularity result we have $\Phi \in \bigoplus_{i=1}^NH^2(\mathbb R^3)$ (cf. \cite{Li}). In order to classify the solutions obtained by actual calculations, it would be important how $A_k(E)$ is arranged in the whole space $\bigoplus_{i=1}^NH^2(\mathbb R^3)$. As for such a problem we have the following corollary. A point $x$ of a real-analytic space $(A,{}_A\mathscr O)$ is called a regular point, if for a neighborhood $W\subset A$ of $x$ and a ringed space $(A',{}_{A'}\mathscr O)$ of a real-analytic subset $A'\subset V\subset \mathbb R^{\mu}$ isomorphic to $(W,{}_A\mathscr O|_W)$ the following holds. The point $y\in A'$ corresponding to $x$ has a neighborhood $V'_y\subset V$ such that by a real-analytic coordinate change to the coordinates $(y_1,\dots,y_{\mu})$, $A'\cap V'_y$ is expressed as $A'\cap V'_y=\{(y_1,\dots,y_{\mu}): y_{\nu+1}=f_{\nu+1}(y_1,\dots,y_{\nu}),\dots,y_{\mu}=f_{\mu}(y_1,\dots,y_{\nu})\}$, where $1\leq\nu\leq\mu-1$ and $f_{\nu+1},\dots,f_{\mu}$ are real-analytic functions. We denote the set of all regular points of $A$ by $\mathfrak R(A)$. For a subspace $X$ of $\bigoplus_{i=1}^NH^2(\mathbb R^3)$ we denote by $X^{\perp}$ its orthogonal subspace and by $P_X$ the orthogonal projection onto $X$.

\begin{corollary}\label{finitesubm}
Let $E<J(N-1)$. Then for any $\delta>0$ there exists a finite-dimensional subspace $X_{\delta}$ of $\bigoplus_{i=1}^NH^2(\mathbb R^3)$ such that 
\begin{itemize}
\item[(i)] For any $\Phi\in A(E)$ we have $\lVert P_{X_{\delta}^{\perp}}\Phi\rVert_{\bigoplus_{i=1}^NH^2(\mathbb R^3)}<\delta$.\\
\end{itemize}
Moreover,
\begin{itemize}
\item[(ii)] For any connected open subset $D$ of $\mathfrak R(A(E))$ such that the closure $D^a$ satisfy $D^a\subset\mathfrak R(A(E))$, $X_{\delta}$ can be chosen so that $P_{X_{\delta}}D$ will be a real-analytic submanifold of $X_{\delta}$.
\end{itemize}
\end{corollary}

This corollary means that $D$ can be approximated by a submanifold of some finite-dimensional subspace of $\bigoplus_{i=1}^NH^2(\mathbb R^3)$ as accurately as one likes.
\begin{remark}
Although $D$ is locally homeomorphic to a finite-dimensional Euclidean space, it may not be included in any finite-dimensional subspace of $\bigoplus_{i=1}^NH^2\newline (\mathbb R^3)$. This can be seen from the following proposition whose proof is irrelevant to the rest of the contents of this paper.
\end{remark}

We denote by $l^2$ the set of all sequences $a=(a_k)_{k=1,2,\dots}$ of complex numbers such that $\sum_{k=1}^{\infty}\lvert a_k\rvert^2<\infty$. Then $l^2$ is a Hilbert space with the innerproduct $(a,b)=\sum_{k=1}^{\infty}a_k^*b_k$. For a subset $M$ of $l^2$ we consider the relative topology induced by the topology of $l^2$.
\begin{proposition}
For any $m\in \mathbb N$ there exists a subset $M$ of $l^2$ homeomorphic to an open unit ball in $\mathbb R^m$ and not included in any finite-dimensional subspace of $l^2$. 
\end{proposition}

\begin{proof}
Let $f_k: B \to B,\ k=1,2,\dots$ be continuous mappings, where $B$ is an open unit ball in $\mathbb R^m$. We also assume that $f_1$ is a homeomorphism. Define $M\subset l^2$ by $M:=\{(f_1(x),2^{-1}f_2(x),\dots,2^{-k}f_k(x),\dots)\in l^2: x\in B\},$
and let $f: B\to M$ be the mapping defined by $f(x):=(f_1(x),2^{-1}f_2(x),\dots,2^{-k}f_k(x),\dots).$
Then we can easily see that $f$ is continuous. Since $f_1$ is injective, $f$ is also injective. For any $q\in\mathbb N$ let $P_q:l^2\to\mathbb C^{mq}$ be a mapping defined by $P_q(y_1,y_2,\dots,y_q,y_{q+1},\dots)=(y_1,\dots,y_q),\ y_i\in \mathbb C^{m}$. Then we have $f^{-1}=f^{-1}_1\circ P_1$ on $M$.
Hence the continuity of $f^{-1}$ follows from the continuity of $f_1^{-1}$.
Thus $M$ is homeomorphic to $B$ by the homeomorphism $f^{-1}$.

We shall choose $f_k$ so that $M$ will have infinitely many linearly independent elements in $l^2$. Then $M$ can not be included in any finite-dimensional subspace of $l^2$. For this purpose we have only to show that there exists a sequence $( x_k)_{k=1,2,\dots},\ x_k\in B$ and $f_k,\ k=1,2\dots$ such that for any $s\in\mathbb N$, $P_{s}\circ f(x_1),P_{s}\circ f(x_2),\dots,P_{s}\circ f(x_s)$ are linearly independent in $\mathbb C^{ms}$. We denote the assertion for $s=t$ by $a(t)$. Then $a(1)$ is obvious. Assuming that $a(t)$ holds, we shall choose $f_{t+1}$ and $x_{t+1}$ so  that $a(t+1)$ holds. We choose $f_{t+1}$ and $x_{t+1}\in B$ satisfying $f_{t+1}(x_1)=f_{t+1}(x_2)=\dotsm=f_{t+1}(x_t)=0$ and $f_{t+1}(x_{t+1})\neq 0$. It is obvious that such a continuous mapping $f_{t+1}$ and a point $x_{t+1}\in B$ exist. Then using the assumption $a(t)$ we can easily see that $P_{t+1}\circ f(x_1),\dots,P_{t+1}\circ f(x_{t+1})$ are linearly independent, i.e. $a(t+1)$ holds. This completes the proof.
\end{proof}

\section{Some preliminaries}\label{thirdsec}
\subsection{Real-analytic operators in Banach space}
In this subsection following \cite{FNSS} we introduce the real-analytic operators and their property. Let $X$ and $Y$ be real Banach spaces.  We denote the norm of $X$ by $\lVert \cdot\rVert$.

\begin{definition}
Let $D$ be an open subset of $X$. The mapping $F:D\to Y$ is said to be real-analytic on $D$ if the following conditions are fulfilled:
\begin{itemize}
\item[(i)] For each $x\in D$ there exist Fr\'echet derivatives of arbitrary orders $d^mF(x,\dots)$.
\item[(ii)] For each $x\in D$ there exists $\delta>0$ such that for any $h\in X$ satisfying $\lVert h\rVert<\delta$ one has
\begin{equation}\label{meq3.1}
F(x+h)=\sum_{m=0}^{\infty}\frac{1}{m!}d^mF(x,h^m),
\end{equation}
(the convergence being locally uniform and absolute), where $h^m := [h, \dots ,h]$ ($m$-times).
\end{itemize}
\end{definition}

The composition of two real-analytic operators is real-analytic (cf. \cite{FNSS, FNSS2}).

\begin{lemma}[{\cite[Proposition 2.1]{FNSS} see also \cite[Lemma 3R]{FNSS2}}]\label{Implicit}
Let $X,Y,Z$ be real Banach spaces, $G\subset X\times Y$ an open set and $[x_0,y_0]\in G$. Let $F: G\to Z$ be a real-analytic mapping such that $[F_y'(x_0,y_0)]^{-1}$ exists and $F(x_0,y_0)=0$. (We denote by
$$F'_y(x_0,y_0)h=\lim_{t\to0}[F(x_0,y_0+th)-F(x_0,y_0)]/t,$$
for $h\in Y$ the partial derivative by $y$. Under our assumptions this derivative exists in the Fr\'echet sense.)

Then there exist a neighborhood $U(x_0)$ in $X$ of the point $x_0$ and a neighborhood $U(y_0)$ in $Y$ of the point $y_0$  such that $U(x_0)\times U(y_0)\subset G$ and there exists one and only one mapping $y: U(x_0)\to U(y_0)$ for which $F(x,y(x))=0$ on $U(x_0)$. Moreover, $y$ is a real-analytic mapping on $U(x_0)$.
\end{lemma}

\subsection{Fr\'echet derivatives of an auxiliary functional at critical points}
We introduce an auxiliary functional whose critical points automatically satisfy the normalizing constraints $\langle \varphi_i,\varphi_i\rangle=1,\ i=1,\dots,N$. This functional plays a crucial role in our study of sets of solutions to the Hartree-Fock equation.
Denote by $Y_1:=(\bigoplus _{i=1}^NH^2(\mathbb R^3))\bigoplus \mathbb R^N$ and $Y_2:=(\bigoplus _{i=1}^NL^2(\mathbb R^3))\bigoplus \mathbb R^N$ the direct sums of Banach spaces regarding the sets $H^2(\mathbb R^3)$ and $L^2(\mathbb R^3)$ of complex-valued functions as real Banach spaces with respect to multiplication by real numbers. We define a functional $f:Y_1\to \mathbb R$ by $f(\Phi,\mathbf e):=\mathcal E(\Phi)-\sum_{i=1}^N\epsilon_i(\lVert \varphi_i\rVert^2-1),$
and a bilinear form $\langle\langle\cdot,\cdot\rangle\rangle$ on $Y_1$ and $Y_2$ by
$$\langle\langle[\Phi^1,\mathbf e^1],[\Phi^2,\mathbf e^2]\rangle\rangle:=\sum_{i=1}^N2\mathrm{Re}\, \langle \varphi_i^1,\varphi_i^2\rangle+\sum_{i=1}^N\epsilon_i^1\epsilon_i^2,$$
where $[\Phi^j,\mathbf e^j] \in Y_j,\ j=1,2$ with $\Phi^j={}^t(\varphi^j_1,\dots,\varphi^j_N)$ and $\mathbf e^j=(\epsilon^j_1,\dots,\epsilon_N^j)$.
We also define a mapping $F: Y_1\to Y_2$ by
$$F(\Phi,\mathbf e):=[{}^t(F_1(\Phi,\mathbf e),\dots,F_N(\Phi,\mathbf e)),(1 - \lVert \varphi_1\rVert^2,\dots,1 - \lVert \varphi_N\rVert^2)],$$
where $F_i(\Phi,\mathbf e):=\mathcal F(\Phi)\varphi_i-\epsilon_i\varphi_i$.
Then we have
$df([\Phi,\mathbf e],[\tilde \Phi,\tilde{\mathbf e}])=\langle\langle[\tilde \Phi,\tilde{\mathbf e}],F(\Phi,\mathbf e)\rangle\rangle$. Moreover, by \eqref{myeq1.3} we can see that if $\Phi={}^t(\varphi_1,\dots,\varphi_N)$ is a solution to the Hartree-Fock equation associated with an orbital energy $\mathbf e$ satisfying the constraints $\langle \varphi_i,\varphi_i\rangle =1,\ 1\leq i\leq N$, we have $F(\Phi,\mathbf e)=0$, and hence $[\Phi,\mathbf e]$ is a critical point of $f(\Phi,\mathbf e)$. Conversely, all critical points of $f(\Phi,\mathbf e)$ are pairs of such a solution and an orbital energy.
We can also see that $F$ is real-analytic as follows. All terms in  $F(\Phi,\mathbf e)$ contain up to three factors such as $\varphi_i$ and $\varphi_i^*$, and they contain up to one factor like $\epsilon_i$. Moreover, there is no term containing more than four factors of $\varphi_i$, $\varphi_i^*$ and $\epsilon_i$. Thus Fr\'echet derivatives of $F(\Phi,\mathbf e)$ higher than third order vanish, and the right-hand side of \eqref{meq3.1} is a finite sum up to third order terms. Hence the convergence is obviously locally uniform and absolute. Set $x=[\Phi,\mathbf e]$ and $h=[\tilde \Phi,\tilde {\mathbf e}]$. $k$-th order differentiation of the $k$-th order terms of $F(x+h)$ with respect to the components of $h$ and substitution of $h^k$ into the derivative reproduce that term of $F(x+h)$ multiplied by $k!$, where $k!$ comes from the number of orders of differentiation of factors (note that this calculation is the same as that for monomials of several variables). Thus the right-hand side of \eqref{meq3.1} equals $F(x+h)$, which means that $F(\Phi,\mathbf e)$ is real-analytic.
Let us denote by $F'(\Phi,\mathbf e)$ the Fr\'echet derivative of $F$ at $[\Phi,\mathbf e]$.

\begin{lemma}\label{LpM}
Let $\Phi^0={}^t(\varphi^0_1\dots,\varphi^0_N)$ be a solution to the Hartree-Fock equation \eqref{myeq1.3} with the constraints $\langle \varphi_i,\varphi_i\rangle=1,\ 1\leq i\leq N$ and not necessarily satisfying the constraints $\langle \varphi_i,\varphi_j\rangle=0,\ i\neq j$ and $\mathbf e^0=(\epsilon^0_1,\dots,\epsilon^0_N)$ be the associated orbital energy such that $\epsilon_i^0\leq-\epsilon,\ 1\leq i\leq N$ for some $\epsilon>0$. Then we have $F'(\Phi^0,\mathbf e^0)=L+M,$
where $L$ is an isomorphism of $Y_1$ onto $Y_2$ and $M$ is a compact operator.
\end{lemma}
Lemma \ref{LpM} has been proved in the proof of \cite[Theorem 2.1]{As} for the Hartree-Fock equation with the constraints $\langle \varphi_i,\varphi_j\rangle=\delta_{ij},\ 1\leq i, j\leq N$. However, we have not used the orthogonal constraints $\langle \varphi_i,\varphi_j\rangle=0,\ i\neq j$ at all in the proof. We give the proof of Lemma \ref{LpM} in an isolated form in the appendix for reader's convenience.

\subsection{Compactness of sets of solutions to the Hartree-Fock equation}
In this subsection we consider the compactness of sets of all solutions to the Hartree-Fock equation with the constraints $\langle \varphi_i,\varphi_i\rangle=1,\ 1\leq i\leq N$ the orbital energies of which satisfy the condition $\epsilon_i\leq-\epsilon,\ 1\leq i\leq N$ for some $\epsilon>0$.
\begin{lemma}[{\cite[Lemma 3.3]{As}}]\label{convseq}
Let $\mathbf e^m:=(\epsilon_1^m,\dots,\epsilon_N^m),\ m=1,2,\dots$ be a sequence of orbital energies converging to $\mathbf e^{\infty}:=(\epsilon_1^{\infty},\dots,\epsilon_N^{\infty})\in(-\infty,0)^N$ and $\Phi^m := {}^t(\varphi_1^m,\dots,\varphi_N^m)$ be the associated solutions to the Hartree-Fock equation \eqref{myeq1.3} with the constraints $\langle \varphi_i,\varphi_i\rangle=1,\ 1\leq i\leq N$ and not necessarily satisfying the constraints $\langle \varphi_i,\varphi_j\rangle=0,\ i\neq j$. Then $\mathbf e^{\infty}$ is an orbital energy and there exists a subsequence of $\Phi^m$ converging in $\bigoplus_{i=1}^NH^2(\mathbb R^3)$ to a solution of the Hartree-Fock equation with the constraints $\langle \varphi_i,\varphi_i\rangle=1,\ 1\leq i\leq N$ associated with $\mathbf e^{\infty}$.
\end{lemma}
\begin{remark}
Lemma \ref{convseq} for the Hartree-Fock equation with the constraints $\langle \varphi_i,\varphi_j\rangle \newline=\delta_{ij}$ has been proved  in \cite[Lemma 3.3]{As} using the Rellich selection theorem and uniform exponential decay of a sequence of solutions. However, the proof has not used the orthogonal constraints $\langle \varphi_i,\varphi_j\rangle=0,\ i\neq j$ at all. Thus we can see that the result holds also for the Hartree-Fock equation with only the normalizing constraints $\langle \varphi_i,\varphi_i\rangle=1,\ 1\leq i\leq N$ as in Lemma \ref{convseq}.
\end{remark}

By Lemma \ref{convseq} we can see that the compactness of the set of solutions follows from boundedness and a uniform negative upper bound of the set of orbital energies. A lower bound for the orbital energies is obtained by the following lemma.
\begin{lemma}\label{olb}
Any orbital energy $\mathbf e=(\epsilon_1,\dots,\epsilon_N)$ of the Hartree-Fock equation \eqref{myeq1.3} with the constraints $\langle \varphi_i,\varphi_i\rangle=1,\ 1\leq i\leq N$ and not necessarily satisfying the constraints $\langle \varphi_i,\varphi_j\rangle=0,\ i\neq j$ satisfies $\epsilon_i\geq \inf\sigma(h)>-\infty,\ 1\leq i\leq N$, where $\sigma(h)$ is the spectrum of $h$.
\end{lemma}

\begin{proof}
First, for $\Phi={}^t(\varphi_1,\dots,\varphi_N)\in \bigoplus_{i=1}^NH^2(\mathbb R^3)$ we shall prove
\begin{equation}\label{meq4.1}
R^{\Phi}-S^{\Phi}\geq0,
\end{equation}
in the sense of inequality between operators.
Since we can write $R^{\Phi}-S^{\Phi} =\sum_{i=1}^N(Q_{ii}^{\Phi}-S_{ii}^{\Phi})$, we only need to prove $Q_{ii}^{\Phi}-S_{ii}^{\Phi} \geq 0,\ 1 \leq i \leq N.$
Let $w \in L^2(\mathbb R^3)$. We define $\hat \Psi_i := 2^{-1/2}(w(x)\varphi_i(y)-\varphi_i(x)w(y)).$
Then it is easily seen that $\langle w ,(Q_{ii}^{\Phi}-S_{ii}^{\Phi})w\rangle = \int \frac{1}{\lvert x-y\rvert}\lvert \hat \Psi_i\rvert^2 dx dy \geq 0.$
Hence we have $Q_{ii}^{\Phi}-S_{ii}^{\Phi} \geq 0$.
Multiplying $\varphi_i^*$ to the Hartree-Fock equation \eqref{myeq1.3} and integrating the both sides we obtain by \eqref{meq4.1}
$$\epsilon_i=\langle\varphi_i,h\varphi_i\rangle+\langle\varphi_i,R^{\Phi}\varphi_i\rangle-\langle\varphi_i,S^{\Phi}\varphi_i\rangle \geq \langle\varphi_i,h\varphi_i\rangle \geq \inf\sigma(h).$$
\end{proof}

Now we can prove the compactness. For $E\in\mathbb R$ and $\epsilon>0$ let $\tilde B(\epsilon)$ be the set of all critical points $[\Phi,\mathbf e]$ of the functional $f(\Phi,\mathbf e)$ (i.e. solutions to the Hartree-Fock equation with the constraints $\langle \varphi_i,\varphi_i\rangle=1,\ 1\leq i\leq N$) such that $\epsilon_i\leq-\epsilon$, $1\leq i\leq N$, where $\mathbf e=(\epsilon_1,\dots,\epsilon_N)$.
\begin{lemma}\label{CompactA}
For any $\epsilon>0$, $\tilde B(\epsilon)$ is a compact set in $Y_1$.
\end{lemma}
\begin{proof}
Let $[\Phi^m,\mathbf e^m]$ be a sequence in $\tilde B(\epsilon)$. By Lemma \ref{olb} and $\epsilon_i\leq-\epsilon$, $1\leq i\leq N$ we can see that there exists a converging subsequence $\mathbf e^{m_k}$ of $\mathbf e^m$, and the limit $\mathbf e^{\infty}=(\epsilon^{\infty}_1,\dots,\epsilon^{\infty}_N)$ satisfies $\epsilon^{\infty}_i\leq -\epsilon<0$. Therefore, by Lemma \ref{convseq} there exists a further subsequence of $[\Phi^{m_k},\mathbf e^{m_k}]$ converging to some $[\Phi^{\infty},\mathbf e^{\infty}]\in \tilde B(\epsilon)$, which means that $\tilde B(\epsilon)$ is a compact set.
\end{proof}

\section{Proof of the main results}\label{fourthsec}
For the proof of the main results we use the notion of semi-analytic subsets and the property of local finiteness of its connected components (cf. \L ojasiewicz \cite{Lo,Lo2}). Let $M$ be a real-analytic manifold. We say that a subset $A\subset M$ is semi-analytic if for any point $x\in M$ there exists a neighborhood $U\subset M$ such that $A\cap U$ is determined by a system of a finite number of inequalities having the form $f> 0$ or $f\geq0$, where by $f$ we mean functions real-analytic in $U$.

\begin{proof}[Proof of Theorem \ref{asub}]
\noindent\textit{Step 1.} Fix some $E<J(N-1)$. First, note that since the orbital energy $\mathbf e=(\epsilon_1,\dots,\epsilon_N)$ is determined by
\begin{equation}\label{myequ4.1}
\epsilon_i=\langle \varphi_i,\mathcal F(\Phi)\varphi_i\rangle,
\end{equation}
from $\Phi$, we can identify $A(E)$ with $\{[\Phi,\mathbf e]: \epsilon_i=\langle \varphi_i,\mathcal F(\Phi)\varphi_i\rangle,\ \Phi\in A(E)\}\subset Y_1$.
Let $z^0=(\Phi^0,\mathbf e^0)\in A(E)$. Then clearly we have $F(\Phi^0,\mathbf e^0)=0$ and by Koopmans' theorem we have $\epsilon^0_i\leq E-J(N-1)<0,\ 1\leq i\leq N$, where $\mathbf e^0=(\epsilon_1^0,\dots,\epsilon_N^0)$ (see \cite[Proof of Theorem 2.3]{As}). Thus by Lemma \ref{LpM} we have $F'(\Phi^0,\mathbf e^0)=L+M,$
where $L$ is an isomorphism of $Y_1$ onto $Y_2$ and $M$ is a compact operator. Notice here that $F'(\Phi^0,\mathbf e^0)$ is bounded by the definition of the Fr\'echet derivative. Thus $L=F'(\Phi^0,\mathbf e^0)-M$ is also bounded. Therefore, by the open mapping theorem it follows that $L^{-1}$ is bounded. Since $L^{-1}M$ and $ML^{-1}$ are compact, by Atkinson's theorem (see e.g. \cite{Ar}) $I+L^{-1}M$ and $I+ML^{-1}$ are Fredholm operators. Hence we can see that $L+M=L(I+L^{-1}M)=(I+ML^{-1})L$ is also Fredholm. Thus $Z_1:=\mathrm{Ker}\, (L+M)$ is finite dimensional, and there exists a closed subspace $Z_2$ of $Y_1$ such that $Y_1=Z_1\bigoplus Z_2$ and $L+M$ is an isomorphism of $Z_2$ onto the closed subspace $Y_3:=(L+M)(Y_1)$ of $Y_2$. If $Z_1=\{0\}$ (i.e. $Z_2=Y_1$), by the inverse mapping theorem (cf. \cite[Theorem 4.F]{Ze}) $F$ is injective on a neighborhood of $z_0$, and thus $z_0$ is an isolated point of $A(E)$. However, setting $\Phi_c:={}^t(c\varphi_1,\dots,c\varphi_N)$ with $c\in \mathbb C$ such that $\lvert c\rvert=1$, it is easily seen that $[\Phi_c,\mathbf e^0]\in A(E)$. Since $[\Phi_c,\mathbf e^0]\to z^0$ as $c\to1$, $z^0$ can not be an isolated point. Thus it follows that $Z_1\neq\{0\}$. Let us write $z^0=[z_1^0,z_2^0],\ z_i^0\in Z_i \ (i=1,2)$ correspondingly to the decomposition $Y_1=Z_1\bigoplus Z_2$.  Let $P_{Y_3}$ be a projection on $Y_3$. Applying Lemma \ref{Implicit} to $P_{Y_3}\circ F$ we can see that there exists a neighborhood $U_1$ of $z_1^0$, a neighborhood $U_2$ of $z_2^0$ and a real-analytic function $\omega: U_1\to U_2$ such that $z=[z_1,z_2]\in U_1\times U_2$ satisfies $P_{Y_3}\circ F([z_1,z_2])=0$ if and only if $z_2=\omega(z_1)$. Let us denote the projection from $Y_1=Z_1\bigoplus Z_2$ onto $Z_1$ by $P_{Z_1}$. Then by a general argument $P_{Z_1}$ is continuous (cf. \cite[page 167]{Ka}). Remembering that $Z_1$ is a finite-dimensional linear space, we have an isomorphism $v:Z_1\to\mathbb R^{\mathrm{dim}\, Z_1}$ which maps $z_1$ to the coefficients in the expression of $z_1$ as a linear combination of a basis in $Z_1$. Set $\mathcal A:=\{z\in Y_1: P_{Y_3}\circ F(z)=0\}$. Let us define a map $u:\mathcal A\cap(U_1\times U_2)\to \mathbb R^{\mathrm{dim}\, Z_1}$ by $u(z)=v\circ {}P_{Z_1}(z)$ for $z\in\mathcal A\cap(U_1\times U_2)$. and define an open subset $V$ of $\mathbb R^{\mathrm{dim}\, Z_1}$ by $V:=u(\mathcal A\cap(U_1\times U_2))=v(U_1)$:

\begin{center}
\begin{tikzpicture}[auto]
\node (b) at (0, 2) {$\mathcal A\cap(U_1\times U_2)$}; \node (a) at (1.8, 2) {$\subset Y_1$}; 
\node (c) at (-0.795, 0) {$Z_1\supset$};   \node (d) at (0.005, 0) {$U_1$};
\node (f) at (2.2, 0) {$V$}; \node (e) at (3.3, 0) {$\subset\mathbb R^{\mathrm{dim}\, Z_1}$};
\draw[->] (b) to node {$\scriptstyle {P_{Z_1}}$} (d);
\draw[->] (d) to node {$\scriptstyle v$} (f);
\draw[->] (b) to node {$\scriptstyle u$} (f);
\end{tikzpicture}
\end{center}
Then $u$ is continuous and its inverse given by $u^{-1}(c)=[v^{-1}(c),\omega\circ v^{-1}(c)],\ c\in V$ is also continuous.
Thus $u$ is an homeomorphism from $\mathcal A\cap(U_1\times U_2)$ to $V$.  We regard $u$ as a coordinate function of $\mathcal A\cap(U_1\times U_2)$. Let $\zeta_1,\dots,\zeta_p$ be a basis of $N(F'(\Phi^0,\mathbf e^0)^*)$ which is finite-dimensional because $F'(\Phi^0,\mathbf e^0)$ is Fredholm, where $N(T)$ is the null space of $T$. Then we can see that $A(E)\cap(U_1\times U_2)$ is the subset of $\mathcal A\cap(U_1\times U_2)$ given by $\zeta_j(F(z))=0,\ j=1,\dots,p$, $\langle \varphi_i,\varphi_j\rangle=0,\ i\neq j$ and $\mathcal E(\Phi)=E$ (Note that $\zeta_j(F(z))=0,\ j=1,\dots,p$ implies $F(z)=0$ for $z\in\mathcal A$, since $P_{Y_3}\circ F(z)=P_{Y_3}F(z)=0$ for $z\in\mathcal A$). Since $F$ and $\mathcal E(\Phi)$ are real-analytic and $\varphi_i$ and $\epsilon_i$ depend real-analytically on the local coordinates, $\zeta_j(F(z))$, $\langle \varphi_i,\varphi_j\rangle$ and $\mathcal E(\Phi)$ are real-analytic functions with respect to the coordinates given by $u$. Therefore, $A':=u(A(E)\cap(U_1\times U_2))$ is a real-analytic subset of $V$.
\smallskip

\noindent\textit{Step 2.} By the argument above, $A(E)$ is identified with an analytic subset in a neighborhood of each point $x\in A(E)$. We shall consider the connection of the analytic subsets in an intersection of two neighborhoods.
Let $A(E)\cap(U_1\times U_2)$ and $A(E)\cap(\tilde U_1\times \tilde U_2)$ be neighborhoods as above and have a non-empty intersection. Hereafter, in this proof we will denote the set $V$ corresponding to $A(E)\cap(\tilde U_1\times \tilde U_2)$ by $\tilde V$ and the mapping $u$ corresponding to $A(E)\cap(\tilde U_1\times \tilde U_2)$ by $\tilde u$ etc. Let $\mathscr O$ (resp., $\tilde {\mathscr O}$) be the sheaf of germs of real-analytic functions on $V$ (resp., $\tilde V$). Moreover, let $\mathscr I\subset \mathscr O$ (resp., $\tilde {\mathscr I}\subset \tilde {\mathscr O}$) be the sheaf of ideals of the real-analytic subset $A'$ (resp., $\tilde{A'}$), that is, the sheaf of germs in $\mathscr O$ (resp., $\tilde {\mathscr O}$) that vanish on $A'$ (resp., $\tilde{A'}$). Then $A'$ (resp., $\tilde{A'}$) is a ringed space with the sheaf ${}_{A'}\mathscr O:=(\mathscr O/\mathscr I)|_{A'}$ (resp., ${}_{\tilde{A'}}\mathscr O:=(\tilde{\mathscr O}/\tilde {\mathscr I})|_{\tilde{A'}}$). Set $\check A:=A(E)\cap(U_1\times U_2)\cap(\tilde U_1\times \tilde U_2)$. In order to show that $A(E)$ has a structure as a real-analytic space we have only to prove that the mapping $\tilde u\circ (u^{-1}|_{u(\check A)})$ is a mapping of ringed spaces from $u(\check A)\subset A'$ to $\tilde u(\check A)\subset\tilde{A'}$ (cf. \cite[Chapter V, A, Proposition 7]{GR}),
that is, for any $y\in u(\check A)$ and $\mathbf g\in{}_{\tilde{A'}}\mathscr O_{\tilde u\circ u^{-1}(y)}$ we have $(\tilde u\circ(u^{-1}|_{u(\check A)}))^*_y(\mathbf g) \in {}_{A'}\mathscr O_y$, where ${}_{A'}\mathscr O_{y}$ is the stalk of ${}_{A'}\mathscr O$ at $y$ and $(\tilde u\circ(u^{-1}|_{u(\check A)}))^*_y$ is the mapping between the stalks of sheaves of germs of continuous functions induced by the pullback by $\tilde u\circ(u^{-1}|_{u(\check A)})$ at $y$ (note that there exists a natural injection from ${}_{A'}\mathscr O_{y}$ to the stalk at $y$ of the sheaf of germs of continuous functions on $A'$). In other words, we need to prove that if we choose some representative function $g$ of $\mathbf g$ and some small neighborhood $V'_y\subset V$ of $y$, then $g\circ(\tilde u\circ (u^{-1}|_{u(\check A)\cap V'_y}))$ is a restriction to $u(\check A)\cap V'_y$ of some real-analytic function in $V'_y$.

On $u((\mathcal A\cap(U_1\times U_2))\cap(\tilde{\mathcal A}\cap(\tilde U_1\times \tilde U_2)))$ the mapping $\tilde u\circ u^{-1}$ can be defined by $\tilde u\circ u^{-1}=\tilde v\circ\tilde P\circ\Omega\circ v^{-1}$, where $\Omega(z_1):=[z_1,\omega(z_1)],\ z_1\in U_1$. Since all mappings in the last expression are real-analytic (note that linear mappings are real-analytic), considering the domains of the mappings we can see that $\tilde u\circ u^{-1}$ defined on $u((\mathcal A\cap(U_1\times U_2))\cap(\tilde{\mathcal A}\cap(\tilde U_1\times \tilde U_2)))$ is extended to a real-analytic mapping from $V$ into $\tilde V$. We still denote the mapping by $\tilde u\circ u$. Since $g$ is a real-analytic function in a neighborhood in $\tilde V$ of $\tilde u\circ u(y)$, the function $g\circ((\tilde u\circ u)|_{V'_y})$ is a real-analytic function on $V'_y$ sufficiently small. Since $g\circ(\tilde u\circ (u^{-1}|_{u(\check A)\cap V'_y}))=(g\circ((\tilde u\circ u)|_{V'_y}))|_{u(\check A)\cap V'_y}$, we can see that the mapping $\tilde u\circ (u^{-1}|_{u(\check A)})$ is a mapping of ringed spaces. 
\smallskip

\noindent\textit{Step 3.}
Next, we shall prove that the set of all connected components of $A(E)$ is finite. Assume that there exist infinitely many connected components ${A}_1,{A}_2,\dots$. Let us choose a point $x_j$ from each $A_j$. From Lemma \ref{CompactA}, Koopmans' theorem and continuity of the functionals $\mathcal E(\Phi)$ and $\langle \varphi_i,\varphi_j\rangle$ it follows that there exists a subsequence of $\{x_j\}$ converging to some point $x_{\infty}\in A(E)$. Let $\mathcal A$, $U_1\times U_2$ and $V$ be the set in Step 1 at the point $x_{\infty}$ instead of $z^0$. Then $\mathcal A\cap(U_1\times U_2)$ is homeomorphic to $V$. Since $A(E)\cap(U_1\times U_2) \subset\mathcal A\cap(U_1\times U_2)$ is determined by $\zeta_j(F(z))=0,\ j=1,\dots,p$, $\langle \varphi_i,\varphi_j\rangle=0,\ i\neq j$ and $\mathcal E(\Phi)=E$ as in Step 1, $u(A(E)\cap(U_1\times U_2) )$ is a real-analytic set (and in particular semi-analytic set) of $V$. The set of all connected components of a semi-analytic subset is locally finite (cf. \cite{Lo}\cite[p.93]{Lo2}). However, any neighborhood of $u(x_{\infty})$ contains infinitely many points of $\{u(x_j)\}$ that belong to different connected components $u(A_j\cap (U_1\times U_2))$ of $u(A(E)\cap (U_1\times U_2))$, which is a contradiction and proves the finiteness of the number of connected components of $A(E)$. The compactness of $A(E)$ follows from Koopmans' theorem and Lemma \ref{CompactA}, and the absence of isolated points can be proved in the same way as in the proof of $Z_1\neq\{0\}$ in  Step 1, which completes the proof.
\end{proof}

The proof of Theorem \ref{bsub} is exactly the same as that of Theorem \ref{asub} except that we choose the neighborhood $U_1\times U_2$ of a point $z^0\in B(\epsilon)$ small enough so that $\epsilon_i<-\epsilon,\ 1\leq i\leq N$ in $U_1\times U_2$, that we do not need the constraint $\mathcal E(\Phi)=E$ and that we need the constraint $\epsilon_i<-\epsilon,\ 1\leq i\leq N$ of semi-analytic set in Step 3 because the orbital energy of $x_{\infty}$ can be $-\epsilon$.

\begin{proof}[Proof of Corollary \ref{finitesubm}]
(i) We shall prove the result by contradiction. Let us assume there exists $\delta>0$ such that for any finite-dimensional subspace $X$ of $\bigoplus_{i=1}^NH^2(\mathbb R^3)$ there exists $\Phi\in A(E)$ satisfying $\lVert P_{X^{\perp}}\Phi\rVert_{\bigoplus_{i=1}^NH^2(\mathbb R^3)}\geq \delta$, contrary to the statement. Then we can make a sequence $\Phi^j\in A(E)$ which does not include any converging subsequence as follows. We choose arbitrary $\Phi^1\in A(E)$. Assume that we have already chosen $\Phi^1,\dots,\Phi^l$. Let $\Phi^{l+1}\in A(E)$ be an element such that $\lVert P_{\mathcal L(\Phi^1,\dots,\Phi^l)^\perp}\Phi^{l+1}\rVert_{\bigoplus_{i=1}^NH^2(\mathbb R^3)} \geq\delta$, where $\mathcal L(\Phi^1,\dots,\Phi^l)$ is the subspace spanned by $\Phi^1,\dots,\Phi^l$. Such $\Phi^{l+1}$ exists by the assumption. Then this sequence satisfies $\lVert \Phi^i-\Phi^j\rVert_{\bigoplus_{i=1}^NH^2(\mathbb R^3)}\geq\delta,$ for any $i\neq j$. Thus there is no converging subsequence of $\{\Phi^j\}$, which contradicts the compactness of $A(E)$.
\smallskip

\noindent (ii) \textit{Step 1.} First of all let us recall the conditions for a subset of a manifold to be a submanifold. Let $\mathcal M$ be a differentiable manifold and $\mathcal N$ a subset of $\mathcal M$ which has a structure as a differentiable manifold i.e. a system of local coordinates. Then $\mathcal N$ is a submanifold of $\mathcal M$ if and only if the injection mapping $\iota:\mathcal N\to\mathcal M$ is a homeomorphism onto its image, differentiable, and the differential $\iota_*$ of $\iota$ is injective at all points of $\mathcal N$ (cf. \cite[Definition 2.5]{St2}). In the present case $X_{\delta}$ and $P_{X_{\delta}}D$ correspond to $\mathcal M$ and $\mathcal N$ respectively.

Note that by Lemma \ref{CompactA}, $D$ is a relatively compact subset of $\mathfrak R(A(E))$. Denote by $\hat A$ the connected component of $\mathfrak R(A(E))$ including $D$. Then $\hat A$ is a manifold as in the proof of Theorem \ref{asub}. For a while let us consider the projection of $\hat A$ on a subspace $X_{\delta}$ of $\bigoplus_{i=1}^NH^2(\mathbb R^3)$ instead of that of $D$. We will later restrict the projection to $D$ to find a finite cover of $P_{X_{\delta}} D$ by coordinate neighborhoods.
Before the proof recall that by \eqref{myequ4.1} we have been identifying a point $\Phi\in\hat A$ with a point $[\Phi,\mathbf e]\in Y_1:=(\bigoplus_{k=1}^NH^2(\mathbb R^3))\times\mathbb R^{3N}$. The correspondence is given also by the projection $P_{\Phi}:Y_1\to \bigoplus_{k=1}^NH^2(\mathbb R^3)$ defined by $P_{\Phi}[\Phi,\mathbf e]=\Phi$. We use this identification to avoid the complex notation such as $P_{X_{\delta}}P_{\Phi}\hat A$ when we regard $\hat A$ as a subset of $Y_1$. As a subset of $Y_1$, $\hat A$ has at each point $z^0\in\hat A$ a coordinate neighborhood $\hat A\cap(U_1\times U_2)$ with the mapping $u: \hat A\cap(U_1\times U_2) \to V\subset \mathbb R^{\mu}$ of the local coordinate system such that $u^{-1}$ is real-analytic as in the proof of Theorem \ref{asub}, where $\mu$ is the dimension of $\hat A$. We define $w:V\to P_{X_{\delta}}\hat A\subset X_{\delta}$ by $w:=P_{X_{\delta}}\circ P_{\Phi}\circ u^{-1}$:

\begin{center}
\begin{tikzpicture}[auto]
\node (a) at (-0.45, 3.03) {$\mathbb R^{\mu}\supset V\ni$}; \node (b) at (1, 3) {$u(z)$};  \node (c) at (3.6, 3) {$z=[\Phi,\mathbf e]$}; \node (d) at (5.4, 3.05) {$\in \hat A\subset Y_1$};
\node (e) at (3.6, 1.5) {$\Phi$}; \node (h) at (4.5, 1.5) {$\in P_{\Phi}\hat A$};  \node (i) at (7, 1.5) {(identified with $\hat A$)};
\node (f) at (3.6, 0) {$w(u(z))$}; \node (g) at (5.6, 0) {$\in P_{X_{\delta}}\hat A\subset X_{\delta}$};
\draw[|->] (c) to node {$\scriptstyle  u$} (b);
\draw[|->] (c) to node {$\scriptstyle  P_{\Phi}$} (e);
\draw[|->] (e) to node {$\scriptstyle P_{X_{\delta}}$} (f);
\draw[|->] (b) to node {$\scriptstyle w$} (f);
\end{tikzpicture}
\end{center}
Then we have only to prove the following assertion (A).\\
\begin{itemize}
\item[(A):] For any point $z^0\in\hat A\cap(U_1\times U_2)$ we can make the differential $w_*$ of $w$ at $u(z^0)$ injective, if we add a finite-dimensional subspace of $\bigoplus_{i=1}^NH^2(\mathbb R^3)$ to $X_{\delta}$.\\
\end{itemize}
We show that the result follows from the assertion (A). Notice that $w$ is a real-analytic mapping because the projections $P_{X_{\delta}}$ and $P_{\Phi}$ are real-analytic mappings by the definition. Since $X_{\delta}$ is finite-dimensional, it has a coordinate system $(\xi_1,\dots,\xi_{\nu})$ where $\nu:=\mathrm{dim}\, X_{\delta}$. By the assertion (A), $\mu\leq \nu$ and we may assume that the Jacobian of $\eta:=(\xi_1\circ w,\dots,\xi_{\mu}\circ w)$ is not zero at $u(z^0)$. Let $\hat Z$ be the subspace of $X_{\delta}$ corresponding to the coordinates $(\xi_1,\dots,\xi_{\mu})$ and $P_{\hat Z}$ be the projection onto $\hat Z$. We identify $\hat Z$ with $\mathbb R^{\mu}$. By the real-analytic version of the inverse function theorem, there exists a homeomorphism $\theta$ from a neighborhood $\hat V\subset \hat Z$ of $\eta(z^0)$ to a neighborhood $\tilde V\subset V$ of $u(z^0)$, and both $\theta$ and $\theta^{-1}$ are real-analytic:

\begin{center}
\begin{tikzpicture}[auto]
\node (a) at (0.6, 3.03) {$\tilde V$}; \node (k) at (1.1, 3) {$\subset$}; \node (b) at (1.6, 3) {$V$}; 
\node (e) at (0.4, 1.53) {$\hat U$}; \node (l) at (0.9, 1.5) {$\subset$}; \node (h) at (1.7, 1.5) {$P_{X_{\delta}}\hat A$};  
\node (f) at (0.6, 0.03) {$\hat V$}; \node (m) at (1.1, -0.03) {$\subset$}; \node (g) at (1.6, 0.03) {$\hat Z$};
\draw[->] (a) to node[swap] {$\scriptstyle  w$} (e);
\draw[->] (e) to node[swap] {$\scriptstyle  \hat u$} (f);
\draw[->] (b) to node {$\scriptstyle w$} (h);
\draw[->] (h) to node {$\scriptstyle P_{\hat Z}$} (g);
\draw[->] (f) to[bend left=80] node {$\scriptstyle \theta$} (a);
\draw[->] (b) to[bend left=80] node {$\scriptstyle \eta=P_{\hat Z}\circ w=(\xi_1\circ w,\dots,\xi_{\mu}\circ w)$} (g);
\end{tikzpicture}
\end{center}
We can see that $w\circ\theta$ is real-analytic and since its inverse is given by $P_{\hat Z}$, it is a homeomorphism. Hence we can use $\hat U:=w\circ\theta(\hat V)=w(\tilde V)$ as a coordinate neighborhood of $P_{X_{\delta}}\hat A$ near $w(z^0)$ and $(\xi_1,\dots,\xi_{\mu})$ as a local coordinate system. We denote by $\hat u$ the mapping from $\hat U$ to the local coordinates. In fact  we have $\hat u=P_{\hat Z}|_{\hat U}$. We can choose such a neighborhood at each point of $P_{X_{\delta}}\hat A$ adding a finite-dimensional subspace to $X_{\delta}$. Since $P_{X_{\delta}}D$ is relatively compact in $P_{X_{\delta}}\hat A$, we have a cover of $P_{X_{\delta}} D$ by a finite number of such neighborhoods. Thus we have only to add a finite number of finite-dimensional subspaces to $X_{\delta}$ to obtain such a cover. On the intersections of the coordinate neighborhoods the coordinate transformations are real-analytic by the real-analyticity of $w\circ \theta$ and $P_{\hat Z}$. Thus $P_{X_{\delta}}D$ is a real-analytic manifold.

We have been considering the relative topologies induced by the topology of $\bigoplus_{i=1}^NH^2(\mathbb R^3)$ as topologies of $X_{\delta}$ and $P_{X_{\delta}}D$ so far. Thus the injection mapping $\iota:P_{X_{\delta}}D\to X_{\delta}$ is obviously a homeomorphism onto its image. As we have seen above, each point in $P_{X_{\delta}}D$ has a neighborhood $\hat U$ and a mapping $\hat u$ of the local coordinates on $\hat U$ such that $\iota\circ\hat u^{-1}$ is given by the mapping $(\xi_1,\dots,\xi_{\mu})\mapsto(\xi_1,\dots,\xi_{\nu}),\ \nu\geq \mu$, where $\xi_{\mu+1},\dots,\xi_{\nu}$ are functions of $\xi_1,\dots,\xi_{\mu}$. Thus $\iota_*$ is injective. Hence $P_{X_{\delta}}D$ satisfies all conditions to be a submanifold of $X_{\delta}$.
\smallskip

\noindent\textit{Step 2.} Now it remains to prove (A). Let us write $z^0=[\Phi^0,\mathbf e^0]$. Remember that the local coordinates $u(z)=u([\Phi,\mathbf e])\in \mathbb R^{\mu}$ in $U$ is given by the coefficients of the representation of $P_{Z_1}[\Phi,\mathbf e]$ as the linear combination of some basis in some finite-dimensional subspace $Z_1$ of $Y_1$ as in the proof of Theorem \ref{asub}. Let $c:(-1,1)\to\mathbb R^{\mu}$ be a $C^1$-curve in $\mathbb R^{\mu}$ such that $c(0)=u([\Phi^0,\mathbf e^0])$ and $c'(0)\neq 0$. $c'(0)$ is identified with a tangent vector at $u([\Phi^0,\mathbf e^0])$. Let us define a curve $[\Phi(t),\mathbf e(t)]$ in $Y_1$ by $[\Phi(t),\mathbf e(t)]=(u)^{-1}(c(t))$. If $\Phi'(0)=0$, then by \eqref{myequ4.1} $\mathbf e'(0)=0$ and thus $\frac{d}{dt}\{(u)^{-1}(c(t))\}\big |_{t=0}=0$. However, since $P_{Z_1}\circ (u)^{-1}$ is a restriction to $V\subset \mathbb R^{\mu}$ of an isomorphism between $\mathbb R^{\mu}$ and $Z_1$, this means $c'(0)=0$, which contradicts the assumption $c'(0)\neq 0$. Thus we can see that $\Phi'(0)\neq 0$. If we choose $X_{\delta}$ so that we will have  $\Phi'(0)\in X_{\delta}$ for all $\Phi'(0)$ corresponding to tangent vectors $c'(0)$ at $u(z^0)$, we can see that $w_*(c'(0))=\Phi'(0)$ from the definition $w=P_{X_{\delta}}\circ P_{\Phi}\circ u^{-1}$. Hence $w_*$ is injective, which completes the proof.
\end{proof}

\def\thesection{\Alph{section}}
\setcounter{section}{0}
\section{Appendix}
In this appendix we prove Lemma \ref{LpM}. In the proof it is supposed that $\Phi^0$ and $\mathbf e^0$ satisfy the assumption in Lemma \ref{LpM}.
\begin{proof}[Proof of Lemma \ref{LpM}]
We define $\tilde F: \bigoplus_{i=1}^NH^2(\mathbb R^3)\to\bigoplus_{i=1}^NL^2(\mathbb R^3)$ by
$$\tilde F(\Phi) := {}^t(F_1(\Phi,\mathbf e^0),\dots,F_N(\Phi,\mathbf e^0)).$$
We shall consider the Fr\'echet derivative of $\tilde F(\Phi)$. Let us define $
R^{\Phi}_i(x):=\sum_{j\neq i}\int\lvert x-y\rvert^{-1}\varphi_j^*(y)\varphi_j(y) dy,$
and $S^{\Phi}_i:=\sum_{j\neq i} S_{jj}^{\Phi}.$
Then by $Q_{ii}^{\Phi}\varphi_i=S_{ii}^{\Phi}\varphi_i$, $F_i(\Phi,\mathbf e^0)=\mathcal F(\Phi)\varphi_i-\epsilon^0_i\varphi_i$ is rewritten as
\begin{equation}\label{meqa.1}
F_i(\Phi,\mathbf e^0)=h\varphi_i+R_i^{\Phi}\varphi_i-S_i^{\Phi}\varphi_i-\epsilon_i^0\varphi_i,
\end{equation}
For a mapping $G: \bigoplus_{i=1}^NH^2(\mathbb R^3)\to L^2(\mathbb R^3)$ and $w_j \in H^2(\mathbb R^3)$ we denote by
$$G'_{j}w_j:=\lim_{t\to0}[G(\varphi_1^0,\dots,\varphi_j^0+tw_j,\dots,\varphi_N^0)-G(\varphi_1^0,\dots,\varphi_j^0,\dots,\varphi_N^0)]/t,$$
the partial derivative of $G$ by $w_j$ at $\Phi^0$. Then by direct calculations we have
\begin{equation}\label{meqa.2}
\begin{split}
&[R_i^{\Phi}\varphi_i]'_{i}w_i=R_i^{\Phi^0}w_i,\\
&[R_i^{\Phi}\varphi_i]'_{j}w_j=S_{ij}^{\Phi^0}w_j+\bar S_{ij}^{\Phi^0}w_j,\ j\neq i,
\end{split}
\end{equation}
and
\begin{equation}\label{meqa.3}
\begin{split}
&[S_i^{\Phi}\varphi_i]'_{i}w_i=S_i^{\Phi^0}w_i,\\
&[S_i^{\Phi}\varphi_i]'_{j}w_j=Q_{ij}^{\Phi^0}w_j+\bar S_{ji}^{\Phi^0}w_j,\ j\neq i,
\end{split}
\end{equation}
where $(\bar S_{ij}^{\Phi}w)(x):=\left(\int\lvert x-y\rvert^{-1}w^*(y)\varphi_j(y)dy\right)\varphi_i(x).$

Set $W:={}^t(w_1,\dots,w_N)$. We define operators $\mathcal R,\mathcal Q: \bigoplus_{i=1}^NH^2(\mathbb R^3)\to\bigoplus_{i=1}^NL^2\newline (\mathbb R^3),$
by $(\mathcal RW)_i:=R_i^{\Phi^0}w_i,$ and $(\mathcal QW)_i:=\sum_{j\neq i}Q_{ij}^{\Phi^0}w_j.$
We shall show that $\mathcal R-\mathcal Q$ is a positive definite operator as an operator on the Hilbert space $\bigoplus_{i=1}^NL^2(\mathbb R^3)$. 
We use the notation $[\tilde i\tilde j\vert kl]$ defined by $[\tilde i\tilde j\vert  k l]:=\int \lvert x-y\rvert^{-1}w_i^*(x)w_j(x)(\varphi_k^0)^*(y)\varphi_l^0(y)dxdy.$
Then we can calculate as $\langle W,(\mathcal R-\mathcal Q)W\rangle=\sum_{i=1}^N\sum_{j\neq i}\{[\tilde i\tilde i\vert jj]-[\tilde i\tilde j\vert ji]\}.$
Moreover, we can calculate as
\begin{align*}
&2^{-1}\sum_{i=1}^N\sum_{j\neq i}\int dx_1dx_2\lvert x_1-x_2\rvert^{-1}\lvert w_i(x_1)\varphi_j^0(x_2)-w_j(x_1)\varphi_i^0(x_2)\rvert^2\\
&\quad=2^{-1}\sum_{i=1}^N\sum_{j\neq i}\{[\tilde i\tilde i\vert jj]+[\tilde j\tilde j\vert ii]-[\tilde i\tilde j\vert ji]-[\tilde j\tilde i\vert ij]\}\\
&\quad=\sum_{i=1}^N\sum_{j\neq i}\{[\tilde i\tilde i\vert jj]-[\tilde i\tilde j\vert ji]\}=\langle W,(\mathcal R-\mathcal Q)W\rangle.
\end{align*}
Since the left-hand side is positive, we can see that $\langle W,(\mathcal R-\mathcal Q)W\rangle\geq0$.

The Fr\'echet derivative of the terms ${}^t(h\varphi_1-\epsilon_1^0\varphi_1,\dots,h\varphi_N-\epsilon_1^0\varphi_N)$ of $\tilde F(\Phi)$ is the operator $\mathcal H: \bigoplus_{i=1}^NH^2(\mathbb R^3)\to\bigoplus_{i=1}^NL^2(\mathbb R^3)$ defined by $(\mathcal HW)_i:=(h-\epsilon_i^0)w_i.$
We can write $\mathcal H$ as a matrix of operators as follows: $\mathcal H=\mathrm{diag}\, (h-\epsilon_1^0 ,\dots ,h-\epsilon_N^0),$
where $\mathrm{diag}\, (A_1 ,\dots ,A_N)$ is the diagonal matrix whose diagonal elements are $A_1 ,\dots ,A_N$.
Denote the resolution of identity of $h$ by $E(\lambda)$. Then we can decompose $h$ as $h=hE(-\epsilon/2)+h(1-E(-\epsilon/2)).$
Since $\inf \sigma_{ess}(h)=0$, $hE(-\epsilon/2)$ is a compact operator, where $\sigma_{ess}(h)$ is the essential spectra of $h$. Moreover, we have an inequality $h(1-E(-\epsilon/2))\geq -\epsilon/2$ of operators. Since $\epsilon_i^0\leq-\epsilon$, we obtain $h(1-E(-\epsilon/2))-\epsilon_i^0 \geq \epsilon/2$. Therefore, $\mathcal H$
is decomposed as a sum
\begin{equation}\label{meqa.4}
\mathcal H=\mathcal H_1+\mathcal H_2,
\end{equation}
of a positive definite operator $\mathcal H_1 := \mathrm{diag}\, (h(1-E(-\epsilon/2))-\epsilon_1^0 ,\dots ,h(1-E(-\epsilon/2))-\epsilon_N^0) \geq \epsilon/2,$
 and a compact operator $\mathcal H_2 := \mathrm{diag}\, (hE(-\epsilon/2) ,\dots ,hE(-\epsilon/2)).$

If we also define $\mathcal S$ and $\bar{\mathcal S}$ by $(\mathcal SW)_i=\sum_{j\neq i}S_{ij}^{\Phi^0}w_j -S_i^{\Phi^0}w_i$ and $(\bar{\mathcal S}W)_i=\sum_{j\neq i}(\bar S_{ij}^{\Phi^0}w_j-\bar S_{ji}^{\Phi^0}w_j),$
by \eqref{meqa.1}-\eqref{meqa.4} the Fr\'echet derivative of $\tilde F$ at $\Phi^0$ is
$$\tilde F'(\Phi^0)=\mathcal H_1+\mathcal H_2+\mathcal R-\mathcal Q+\mathcal S+\bar{\mathcal S}=\mathcal L+\mathcal M,$$
where $\mathcal L:=\mathcal H_1+\mathcal R-\mathcal Q$ and $\mathcal M:=\mathcal H_2+\mathcal S+\bar{\mathcal S}$.
Since $\mathcal R-\mathcal Q$ is positive definite and $\mathcal H_1 \geq \epsilon/2$, we have $\mathcal L \geq \epsilon/2$ and thus, $\mathcal L$ is a map onto $\bigoplus_{i=1}^NL^2(\mathbb R^3)$ and invertible. Therefore, $\mathcal L$ is an isomorphism. Moreover, since $S_{ij}^{\Phi^0}$ and $\bar S_{ij}^{\Phi^0}$ are integral operators of Hilbert-Schmidt type, they are compact operators. Hence $\mathcal S+\bar{\mathcal S}$ is compact and $\mathcal M$ is also compact.

Setting $\hat F(\mathbf e) := {}^t(F_1(\Phi^0,\mathbf e),\dots,F_N(\Phi^0,\mathbf e)),$
we can see that
\begin{align*}
F'(\Phi^0,\mathbf e^0)[\Phi,\mathbf e]&= [\tilde F'(\Phi^0)\Phi+\hat F'(\mathbf e^0)\mathbf e,- 2\mathrm{Re}\, \langle\varphi_1,\varphi_1^0\rangle,\dots, - 2\mathrm{Re}\, \langle\varphi_N,\varphi_N^0\rangle]\\
&=L[\Phi,\mathbf e]+M[\Phi,\mathbf e],
\end{align*}
where $L[\Phi,\mathbf e]:=[\mathcal L\Phi,\mathbf e]$,
$$M[\Phi,\mathbf e]:=[\mathcal M\Phi-\mathbf e\Phi^0,- 2\mathrm{Re}\, \langle\varphi_1,\varphi_1^0\rangle-\epsilon_1,\dots, - 2\mathrm{Re}\, \langle\varphi_N,\varphi_N^0\rangle-\epsilon_N],$$
and $\mathbf e\Phi^0 := {}^t(\epsilon_1\varphi_1^0,\dots,\epsilon_N\varphi_N^0)$. By the properties of $\mathcal L$ and $\mathcal M$ we can easily see that $L$ is an isomorphism and $M$ is a compact operator, which completes the proof.
\end{proof}

\end{document}